\documentclass[12pt]{article}
\usepackage{latexsym}
\usepackage{amssymb}
\usepackage{graphicx}
\usepackage{cite}
\usepackage{color}
\usepackage{ifpdf}

\newtheorem{Theorem}{Theorem}[part]
\newtheorem{Definition}{Definition}[part]
\newtheorem{Proposition}{Proposition}[part]

\newtheorem{Lemma}{Lemma}[part]
\newtheorem{Corollary}{Corollary}[part]

\newtheorem{Conjecture}{Conjecture}

\def \ep{\hbox{ }\hfill$\Box$}

\def\reff#1{{\rm(\ref{#1})}}

\begin{document}
\title{The Largest Laplacian and Signless Laplacian H-Eigenvalues of a Uniform Hypergraph}

\author{
Shenglong Hu \thanks{Email: Tim.Hu@connect.polyu.hk. Department of
Applied Mathematics, The Hong Kong Polytechnic University, Hung Hom,
Kowloon, Hong Kong.},\hspace{4mm} Liqun Qi \thanks{Email:
maqilq@polyu.edu.hk. Department of Applied Mathematics, The Hong
Kong Polytechnic University, Hung Hom, Kowloon, Hong Kong. This
author's work was supported by the Hong Kong Research Grant
Council (Grant No. PolyU 501909, 502510, 502111 and 501212).},
\hspace{4mm} Jinshan Xie \thanks{Email: jinshan0623@sina.com. School of Mathematics and Computer Science, Longyan University, Longyan, Fujian, China.}}

\date{\today}
\maketitle

\begin{abstract}
In this paper, we show that the largest Laplacian H-eigenvalue of a
$k$-uniform nontrivial hypergraph is strictly larger than the
maximum degree when $k$ is even. A tight lower bound for this
eigenvalue is given. For a connected even-uniform hypergraph, this
lower bound is achieved if and only if it is a hyperstar. However,
when $k$ is odd, in certain cases the largest Laplacian H-eigenvalue
is equal to the maximum degree, which is a tight lower bound. On the
other hand, tight upper and lower bounds for the largest signless
Laplacian H-eigenvalue of a $k$-uniform connected hypergraph are
given. For a connected $k$-uniform hypergraph, the upper
(respectively lower) bound of the largest signless Laplacian
H-eigenvalue is achieved if and only if it is a complete hypergraph
(respectively a hyperstar). The largest Laplacian H-eigenvalue is
always less than or equal to the largest signless Laplacian
H-eigenvalue. When the hypergraph is connected, the equality holds
here if and only if $k$ is even and the hypergraph is odd-bipartite.

\vspace{2mm}
\noindent {\bf Key words:}\hspace{2mm} Tensor, H-eigenvalue, hypergraph, Laplacian, signless Laplacian\vspace{1mm}

\noindent {\bf MSC (2010):}\hspace{2mm}
05C65; 15A18
\end{abstract}

\section{Introduction}
\setcounter{Theorem}{0} \setcounter{Proposition}{0}
\setcounter{Corollary}{0} \setcounter{Lemma}{0}
\setcounter{Definition}{0} \setcounter{Remark}{0}
\setcounter{Conjecture}{0}  \setcounter{Example}{0} \hspace{4mm} In
this paper, we study the largest Laplacian and signless Laplacian
H-eigenvalues of a uniform hypergraph. The largest Laplacian and
signless Laplacian H-eigenvalues refer to respectively the largest
H-eigenvalue of the Laplacian tensor and the largest H-eigenvalue of
the signless Laplacian tensor. This work is motivated by the classic
results for graphs \cite{c97,bh11,crs07,zl02,z07}. Please refer to
\cite{l07,hq12a,cd12,q12a,pz12,lqy12,rp09,r09,xc12c,hq13,cpz08,hhlq12,l05,q05,q12b,yy11,xc12a}
for recent developments on spectral hypergraph theory and the
essential tools from spectral theory of nonnegative tensors.

This work is a companion of the recent study on the eigenvectors of the zero Laplacian and signless Laplacian eigenvalues of a uniform hypergraph by Hu and Qi \cite{hq13b}. For the literature on the Laplacian-type tensors for a uniform hypergraph, which becomes an active research frontier in spectral hypergraph theory, please refer to \cite{hq12a,lqy12,xc12c,q12a,hq13,xc12a,hq13b} and references therein.
Among others, Qi \cite{q12a}, and Hu and Qi \cite{hq13} respectively systematically studied the Laplacian and signless Laplacian tensors, and the Laplacian of a uniform hypergraph. These three notions of Laplacian-type tensors are more natural and simpler than those in the literature.

The rest of this paper is organized as follows. Some definitions on
eigenvalues of tensors and uniform hypergraphs are presented in the
next section. The class of hyperstars is introduced. We discuss in
Section 3 the largest Laplacian H-eigenvalue of a $k$-uniform
hypergraph. We show that when $k$ is even, the largest Laplacian
H-eigenvalue has a tight lower bound that is strictly larger than the maximum degree.
Extreme hypergraphs in this case are characterized, which are the
hyperstars. When $k$ is odd, a tight lower bound is exactly the
maximum degree. However, we are not able to characterize the extreme
hypergraphs in this case. Then we discuss the largest signless
Laplacian H-eigenvalue in Section 4. Tight lower and upper bounds
for the largest signless Laplacian H-eigenvalue of a connected
hypergraph are given. Extreme hypergraphs are characterized as well.
For the lower bound, the extreme hypergraphs are hyperstars; and for
the upper bound, the extreme hypergraphs are complete hypergraphs.
The relationship between the largest Laplacian H-eigenvalue and the
largest signless Laplacian H-eigenvalue is discussed in Section 5.
The largest Laplacian H-eigenvalue is always less than or equal to
the largest signless Laplacian H-eigenvalue. When the hypergraph is
connected, the equality holds here if and only if $k$ is even and
the hypergraph is odd-bipartite. This result can help to find the largest Laplacian H-eigenvalue of an even-uniform hypercycle.
Some final remarks are made in the
last section.

\section{Preliminaries}\label{sec-p}
\setcounter{Theorem}{0} \setcounter{Proposition}{0}
\setcounter{Corollary}{0} \setcounter{Lemma}{0}
\setcounter{Definition}{0} \setcounter{Remark}{0}
\setcounter{Conjecture}{0}  \setcounter{Example}{0}
\hspace{4mm} Some definitions of eigenvalues of tensors and uniform hypergraphs are presented in this section.

\subsection{Eigenvalues of Tensors}\label{s-et}
In this subsection, some basic definitions on eigenvalues of tensors are reviewed. For comprehensive references, see \cite{q05,hhlq12} and references therein. Especially, for spectral hypergraph theory oriented facts on eigenvalues of tensors, please see \cite{q12a,hq13}.

Let $\mathbb R$ be the field of real numbers and $\mathbb R^n$ the
$n$-dimensional real space. $\mathbb R^n_+$ denotes the nonnegative
orthant of $\mathbb R^n$. For integers $k\geq 3$ and $n\geq 2$, a
real tensor $\mathcal T=(t_{i_1\ldots i_k})$ of order $k$ and
dimension $n$ refers to a multiway array (also called hypermatrix)
with entries $t_{i_1\ldots i_k}$ such that $t_{i_1\ldots
i_k}\in\mathbb{R}$ for all $i_j\in[n]:=\{1,\ldots,n\}$ and
$j\in[k]$. Tensors are always referred to $k$-th order real tensors
in this paper, and the dimensions will be clear from the content.
Given a vector $\mathbf{x}\in \mathbb{R}^{n}$, ${\cal
T}\mathbf{x}^{k-1}$ is defined as an $n$-dimensional vector such
that its $i$-th element being
$\sum\limits_{i_2,\ldots,i_k\in[n]}t_{ii_2\ldots i_k}x_{i_2}\cdots
x_{i_k}$ for all $i\in[n]$. Let ${\cal I}$ be the identity tensor of
appropriate dimension, e.g., $i_{i_1\ldots i_k}=1$ if and only if
$i_1=\cdots=i_k\in [n]$, and zero otherwise when the dimension is
$n$. The following definition was introduced by Qi \cite{q05}.
\begin{Definition}\label{def-00}
Let $\mathcal T$ be a $k$-th order $n$-dimensional real tensor. For
some $\lambda\in\mathbb{R}$, if polynomial system $\left(\lambda
{\cal I}-{\cal T}\right)\mathbf{x}^{k-1}=0$ has a solution
$\mathbf{x}\in\mathbb{R}^n\setminus\{0\}$, then $\lambda$ is called
an H-eigenvalue and $\mathbf x$ an
H-eigenvector.
\end{Definition}
It is seen that H-eigenvalues are real numbers \cite{q05}. By
\cite{hhlq12,q05}, we have that the number of H-eigenvalues of a
real tensor is finite. By \cite{q12a}, we have that all the tensors
considered in this paper have at least one H-eigenvalue. Hence, we
can denote by $\lambda(\mathcal T)$ (respectively $\mu(\mathcal T)$)
as the largest (respectively smallest) H-eigenvalue of a real tensor
$\mathcal T$.

For a subset $S\subseteq [n]$, we denoted by $|S|$ its cardinality, and $\mbox{sup}(\mathbf x):=\{i\in[n]\;|\;x_i\neq 0\}$ its {\em support}.

\subsection{Uniform Hypergraphs}

In this subsection, we present some essential concepts of uniform hypergraphs which will be used in the sequel. Please refer to \cite{b73,c97,bh11,hq13,q12a} for comprehensive references.

In this paper, unless stated otherwise, a hypergraph means an undirected simple $k$-uniform hypergraph $G$ with vertex set $V$, which is labeled as $[n]=\{1,\ldots,n\}$, and edge set $E$.
By $k$-uniformity, we mean that for every edge $e\in E$, the cardinality $|e|$ of $e$ is equal to $k$. Throughout this paper, $k\geq 3$ and $n\geq k$. Moreover, since the trivial hypergraph (i.e., $E=\emptyset$) is of less interest, we consider only hypergraphs having at least one edge (i.e., nontrivial) in this paper.

For a subset $S\subset [n]$, we denoted by $E_S$ the set of edges $\{e\in E\;|\;S\cap e\neq\emptyset\}$. For a vertex $i\in V$, we simplify $E_{\{i\}}$ as $E_i$. It is the set of edges containing the vertex $i$, i.e., $E_i:=\{e\in E\;|\;i\in e\}$. The cardinality $|E_i|$ of the set $E_i$ is defined as the {\em degree} of the vertex $i$, which is denoted by $d_i$. Two different vertices $i$ and $j$ are {\em connected} to each other (or the pair $i$ and $j$ is connected), if there is a sequence of edges $(e_1,\ldots,e_m)$ such that $i\in e_1$, $j\in e_m$ and $e_r\cap e_{r+1}\neq\emptyset$ for all $r\in[m-1]$. A hypergraph is called {\em connected}, if every pair of different vertices of $G$ is connected. Let $S\subseteq V$, the hypergraph with vertex set $S$ and edge set $\{e\in E\;|\;e\subseteq S\}$ is called the {\em sub-hypergraph} of $G$ induced by $S$. We will denote it by $G_S$. A hypergraph is {\em regular} if $d_1=\cdots=d_n=d$.
A hypergraph $G=(V,E)$ is {\em complete} if $E$ consists of all the possible edges. In this case, $G$ is regular, and moreover $d_1=\cdots=d_n=d={n-1\choose k-1}$.
In the sequel, unless stated otherwise, all the notations introduced above are reserved for the specific meanings.

 For the sake of simplicity, we mainly consider connected hypergraphs in the subsequent analysis. By the techniques in \cite{q12a,hq13}, the conclusions on connected hypergraphs can be easily generalized to general hypergraphs.

The following definition for the Laplacian tensor and signless Laplacian tensor was proposed by Qi \cite{q12a}.
\begin{Definition}\label{def-l}
Let $G=(V,E)$ be a $k$-uniform hypergraph. The {\em adjacency tensor} of $G$ is defined as the $k$-th order $n$-dimensional tensor $\mathcal A$ whose $(i_1 \ldots i_k)$-entry is:
\begin{eqnarray*}
a_{i_1 \ldots i_k}:=\left\{\begin{array}{cl}\frac{1}{(k-1)!}&if\;\{i_1,\ldots,i_k\}\in E,\\0&\mbox{otherwise}.\end{array}\right.
\end{eqnarray*}
Let $\mathcal D$ be a $k$-th order $n$-dimensional diagonal tensor with its diagonal element $d_{i\ldots i}$ being $d_i$, the degree of vertex $i$, for all $i\in [n]$. Then $\mathcal L:=\mathcal D-\mathcal A$ is the {\em Laplacian tensor} of the hypergraph $G$, and $\mathcal Q:=\mathcal D+\mathcal A$ is the {\em signless Laplacian tensor} of the hypergraph $G$.
\end{Definition}

In the following, we introduce the class of hyperstars.
\begin{Definition}\label{def-bi-hm}
Let $G=(V,E)$ be a $k$-uniform hypergraph. If there is a disjoint partition of the vertex set $V$ as $V=V_0\cup V_1\cup\cdots\cup V_d$ such that $|V_0|=1$ and $|V_1|=\cdots=|V_d|=k-1$, and $E=\{V_0\cup V_i\;|\;i\in [d]\}$, then $G$ is called a {\em hyperstar}. The degree $d$ of the vertex in $V_0$, which is called the {\em heart}, is the {\em size} of the hyperstar. The edges of $G$ are {\em leaves}, and the vertices other than the heart are vertices of leaves.
\end{Definition}
It is an obvious fact that, with a possible renumbering of the vertices, all the hyperstars with the same size are identical. Moreover, by Definition \ref{def-00}, we see that the process of renumbering does not change the H-eigenvalues of either the Laplacian tensor or the signless Laplacian tensor of the hyperstar.
The trivial hyperstar is the one edge hypergraph, its spectrum is very clear \cite{cd12}. In the sequel, unless stated otherwise, a hyperstar is referred to a hyperstar having size $d>1$.
For a vertex $i$ other than the heart, the leaf containing $i$ is denoted by $le(i)$. An example of a hyperstar is given in Figure 1.

\begin{figure}[htbp]
\centering
\includegraphics[width=2.0in]{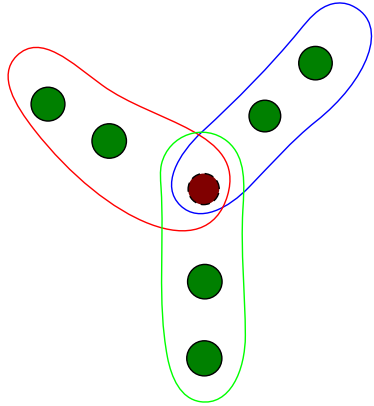}
\caption{An example of a $3$-uniform hyperstar of size $3$. An edge is pictured as a closed curve with the containing solid disks the vertices in that edge. Different edges are in different curves with different colors. The red (also in dashed margin) disk represents the heart.}
\end{figure}

The notions of odd-bipartite and even-bipartite even-uniform hypergraphs are introduced in \cite{hq13b}.
\begin{Definition}\label{def-bi-odd}
Let $k$ be even and $G=(V,E)$ be a $k$-uniform
hypergraph. It is called {\em odd-bipartite} if either it is trivial
(i.e., $E=\emptyset$) or there is a disjoint partition of the vertex
set $V$ as $V=V_1\cup V_2$ such that $V_1,V_2\neq \emptyset$ and
every edge in $E$ intersects $V_1$ with exactly an odd number of
vertices.
\end{Definition}
An example of an odd-bipartite hypergraph is given in Figure 2.

\begin{figure}[htbp]
\centering
\includegraphics[width=2.0in]{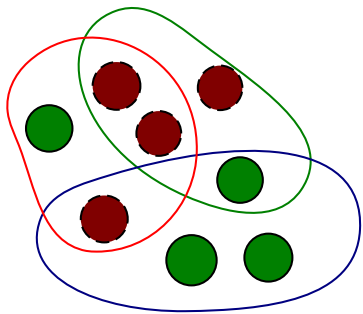}
\caption{An example of an odd-bipartite $4$-uniform hypergraph. The bipartition is clear from the different colors (also the dashed margins from the solid ones) of the disks.}
\end{figure}
\section{The Largest Laplacian H-Eigenvalue}
\setcounter{Theorem}{0} \setcounter{Proposition}{0}
\setcounter{Corollary}{0} \setcounter{Lemma}{0}
\setcounter{Definition}{0} \setcounter{Remark}{0}
\setcounter{Conjecture}{0}  \setcounter{Example}{0} \hspace{4mm}
This section presents some basic facts about the largest Laplacian H-eigenvalue of a uniform hypergraph. We start the discussion on the class of hyperstars.

\subsection{Hyperstars}
Some properties of hyperstars are given in this subsection.

The next proposition is a direct consequence of Definition \ref{def-bi-hm}.
\begin{Proposition}\label{prop-0}
Let $G=(V,E)$ be a hyperstar of size $d>0$. Then except for one vertex $i\in [n]$ with $d_i=d$, we have $d_j=1$ for the others.
\end{Proposition}

By Theorem 4 of \cite{q12a}, we have the following lemma.
\begin{Lemma}\label{lem-00}
Let $G=(V,E)$ be a $k$-uniform hypergraph with its maximum degree
$d>0$ and $\mathcal L=\mathcal D-\mathcal A$ be its Laplacian
tensor. Then $\lambda (\mathcal L)\geq d$.
\end{Lemma}

When $k$ is even and $G$ is a hyperstar, Lemma \ref{lem-00} can be
strengthened as in the next proposition.
\begin{Proposition}\label{prop-1}
Let $k$ be even and $G=(V,E)$ be a hyperstar of size $d>0$ and $\mathcal L=\mathcal D-\mathcal A$ be its Laplacian tensor. Then $\lambda (\mathcal L)>d$.
\end{Proposition}

\noindent {\bf Proof.} Suppose, without loss of generality, that $d_1=d$. Let $\mathbf x\in\mathbb R^n$ be a nonzero vector such that $x_1=\alpha\in\mathbb R$, and $x_2=\cdots=x_n=1$. Then, we see that
\begin{eqnarray*}
\left(\mathcal L\mathbf x^{k-1}\right)_1=d \alpha^{k-1}-d,
\end{eqnarray*}
and for $i\in\{2,\ldots,n\}$
\begin{eqnarray*}
\left(\mathcal L\mathbf x^{k-1}\right)_i=1-\alpha.
\end{eqnarray*}
Thus, if $\mathbf x$ is an H-eigenvector of $\mathcal L$ corresponding to an H-eigenvalue $\lambda$, then we must have
\begin{eqnarray*}
d \alpha^{k-1}-d=\lambda\alpha^{k-1},\;\mbox{and}\; 1-\alpha=\lambda.
\end{eqnarray*}
Hence,
\begin{eqnarray*}
(1-\lambda)^{k-1}(\lambda-d)+d=0.
\end{eqnarray*}
Let $f(\lambda):=(1-\lambda)^{k-1}(\lambda-d)+d$. We have that
\begin{eqnarray*}
f(d)=d>0, \;\mbox{and}\;f(d+1)=(-d)^{k-1}+d<0.
\end{eqnarray*}
Consequently, $f(\lambda)=0$ does have a root in the interval $(d,d+1)$. Hence $\mathcal L$ has an H-eigenvalue $\lambda>d$. The result follows. \ep

The next lemma characterizes H-eigenvectors of the Laplacian tensor of a hyperstar corresponding to an H-eigenvalue which is not one .
\begin{Lemma}\label{lem-0}
Let $G=(V,E)$ be a hyperstar of size $d>0$ and $\mathbf x\in\mathbb R^n$ be an H-eigenvector of the Laplacian tensor of $G$ corresponding to a nonzero H-eigenvalue other than one. If $x_i=0$ for some vertex of a leaf (other than the heart), then $x_j=0$ for all the vertices $j$ in the leaf containing $i$ and other than the heart. Moreover, in this situation, if $h$ is the heart, then $x_h\neq 0$.
\end{Lemma}

\noindent {\bf Proof.} Suppose that the H-eigenvalue is $\lambda\neq 1$. By the definition of eigenvalues, we have that for the vertex $j$
other than the heart and the vertex $i$,
\begin{eqnarray*}
\left(\mathcal L\mathbf x^{k-1}\right)_j=x_j^{k-1}-\prod_{s\in
le(j)\setminus\{j\}} x_s=x_j^{k-1}-0=\lambda x_j^{k-1}.
\end{eqnarray*}
Since $\lambda\neq 1$, we must have that $x_j=0$.

With a similar proof, we get the other conclusion by contradiction, since $h\in le(i)$ for all vertices $i$ of leaves and $\mathbf x\neq 0$. \ep

The next lemma characterizes the H-eigenvectors of the Laplacian tensor of a hyperstar corresponding to the largest Laplacian H-eigenvalue.
\begin{Lemma}\label{lem-1}
Let $G=(V,E)$ be a hyperstar of size $d>1$. Then there is an H-eigenvector $\mathbf z\in\mathbb R^n$ of the Laplacian tensor $\mathcal L$ of $G$ corresponding to $\lambda(\mathcal L)$ satisfying that $|z_i|$ is a constant for $i\in\mbox{sup}(\mathbf z)$ and $i$ being not the heart.
\end{Lemma}

\noindent {\bf Proof.}
Suppose that $\mathbf y\in\mathbb R^{n}$ is an H-eigenvector of $\mathcal L$ corresponding to $\lambda(\mathcal L)$. Without loss of generality, let $1$ be the heart and hence $d_1=d$. Note that, by Lemma \ref{lem-00}, we have that $\lambda(\mathcal L)\geq d>1$.
By Lemma \ref{lem-0}, without loss of generality, we can assume that $\mbox{sup}(\mathbf y)=[n]$ and $y_1>0$. In the following, we construct an H-eigenvector $\mathbf z\in\mathbb R^{n}$ corresponding to $\lambda(\mathcal L)$ from $\mathbf y$ such that
$|z_2|=\cdots=|z_{n}|$.

(I). We first prove that for every leaf $e\in E$, $|y_t|$ is a constant for all $t\in e\setminus\{1\}$.

For an arbitrary but fixed leaf $e\in E$, suppose that $|y_i|=\max\{|y_j|\;|\;j\in e\setminus\{1\}\}$ and $|y_s|=\min\{|y_j|\;|\;j\in e\setminus\{1\}\}$. If $|y_i|=|y_s|$, then we are done. In the following, suppose on the contrary that $|y_i|>|y_s|$. Then, we have
\begin{eqnarray*}
(\lambda(\mathcal L)-1)|y_i|^{k-1}=y_1\prod_{j\in e\setminus\{1,i\}}|y_j|, \;\mbox{and}\;(\lambda(\mathcal L)-1)|y_s|^{k-1}=y_1\prod_{j\in e\setminus\{1,s\}}|y_j|.
\end{eqnarray*}
By the definitions of $|y_i|$ and $|y_s|$, we have $y_1\prod_{j\in e\setminus\{1,i\}}|y_j|< y_1\prod_{j\in e\setminus\{1,s\}}|y_j|$. On the other hand, we have $(\lambda(\mathcal L)-1)|y_i|^{k-1}>(\lambda(\mathcal L)-1)|y_s|^{k-1}$. Hence, a contradiction is derived. Consequently, for every leaf $e\in E$, $|y_t|$ is a constant for all $t\in e\setminus\{1\}$.

(II). We next show that all the numbers in this set
\begin{eqnarray*}
\left\{\alpha_s:=\prod_{j\in e_s\setminus\{1\}}y_j,\;e_s\in
E\right\}
\end{eqnarray*}
are of the same sign.

When $k$ is even, suppose that $y_i<0$ for some $i$. Then
\begin{eqnarray}\label{new-4}
0>(\lambda(\mathcal L)-1)y_i^{k-1}=-y_1\prod_{j\in le(i)\setminus\{1,i\}}y_j.
\end{eqnarray}
Thus, an odd number of vertices in $le(i)$ takes negative values. By \reff{new-4}, we must have that there exists some $i\in e$ such that $y_i<0$ for every $e\in E$. Otherwise, $(\lambda(\mathcal L)-1)y_i^{k-1}>0$, together with $-y_1\prod_{j\in le(i)\setminus\{1,i\}}y_j<0$, would lead to a contradiction.
Hence, all the numbers in this set
\begin{eqnarray*}
\left\{\alpha_s:=\prod_{j\in e_s\setminus\{1\}}y_j,\;e_s\in
E\right\}
\end{eqnarray*}
are negative.

When $k$ is odd, suppose that $y_i<0$ for some $i$. Then
\begin{eqnarray}\label{new-3}
0<(\lambda(\mathcal L)-1)y_i^{k-1}=-y_1\prod_{j\in le(i)\setminus\{1,i\}}y_j.
\end{eqnarray}
Thus, an positive even number of vertices in $le(i)$ takes negative values. Thus, if there is some $s\in le(i)$ such that $y_s>0$, then
\begin{eqnarray*}
0<(\lambda(\mathcal L)-1)y_s^{k-1}=-y_1\prod_{j\in le(s)\setminus\{1,s\}}y_j.
\end{eqnarray*}
Since $s\in le(i)$, we have $le(i)=le(s)$ and $i\in le(s)$. Hence, $y_1\prod_{j\in le(s)\setminus\{1,s\}}y_j>0$. A contradiction is derived. By \reff{new-3}, we must have that there exists some $i\in e$ such that $y_i<0$ for every $e\in E$. Consequently, $y_j<0$ for all $j\neq 1$. Hence, all the numbers in this set
\begin{eqnarray*}
\left\{\alpha_s:=\prod_{j\in e_s\setminus\{1\}}y_j,\;e_s\in
E\right\}
\end{eqnarray*}
are positive.

(III.) We construct the desired vector $\mathbf z$.

If the product $\prod_{j\in e\setminus\{1\}}y_j$ is a constant for every leaf $e\in E$, then take $\mathbf z=\mathbf y$ and we are done.
In the following, suppose on the contrary that the set
\begin{eqnarray*}
\left\{\alpha_s:=\prod_{j\in e_s\setminus\{1\}}y_j,\;e_s\in
E\right\}
\end{eqnarray*}
takes more than one numbers. Let $\mathbf z\in\mathbb R^n$ be the vector such that
\begin{eqnarray*}
z_j=\left(\frac{\sum_{e_s\in E}\alpha_s}{d\alpha_t}\right)^{1 \over
k-1}y_j,\;j\in e_t\setminus\{1\}
\end{eqnarray*}
and $z_1=y_1$. Note that $|z_2|=\cdots=|z_{n_2}|$, since $|y_j|^{k-1}=\alpha_t$ for all $j\in e_t\setminus\{1\}$ and $e_t\in E$. Then
\begin{eqnarray*}
\left(\mathcal L\mathbf z^{k-1}\right)_1&=&dz_1^{k-1}-\sum_{e_s\in E}\prod_{j\in \in e_s\setminus\{1\}}z_j\\
&=&dy_1^{k-1}-\sum_{e_s\in E}\frac{\sum_{e_t\in E}\alpha_t}{d\alpha_s}\prod_{j\in \in e_s\setminus\{1\}}y_j\\
&=&dy_1^{k-1}-\sum_{e_s\in E}\frac{\sum_{e_t\in E}\alpha_t}{d\alpha_s}\alpha_s\\
&=&dy_1^{k-1}- \sum_{e_t\in E}\alpha_t\\
&=&dy_1^{k-1}- \sum_{e_t\in E}\prod_{j\in e_t\setminus\{1\}}y_j\\
&=&\lambda(\mathcal L) y_1^{k-1}\\
&=&\lambda(\mathcal L) z_1^{k-1}.
\end{eqnarray*}
For any $i\neq 1$ with $i\in e_s$ for some $s$, we have
\begin{eqnarray*}
\left(\mathcal L\mathbf z^{k-1}\right)_i&=&z_i^{k-1}-z_1\prod_{j\in e_s\setminus\{1,i\}}z_j\\
&=&\frac{\sum_{e_t\in E}\alpha_t}{d\alpha_s}y_i^{k-1}-y_1\frac{\sum_{e_t\in E}\alpha_t}{d\alpha_s}\prod_{j\in e_s\setminus\{1,i\}}y_j\\
&=&\lambda(\mathcal L)\frac{\sum_{e_t\in E}\alpha_t}{d\alpha_s}y_i^{k-1}\\
&=&\lambda(\mathcal L)z_i^{k-1}.
\end{eqnarray*}
By Definition \ref{def-00}, $\mathbf z$ is an H-eigenvector of $\mathcal L$ corresponding to $\lambda(\mathcal L)$ with the requirement. The result follows. \ep

The next corollary follows directly from the proof of Lemma \ref{lem-1}.
\begin{Corollary}\label{cor-4}
Let $k$ be odd and $G=(V,E)$ be a hyperstar of size $d>1$. If $\mathbf z\in\mathbb R^n$ is an H-eigenvector of the Laplacian tensor $\mathcal L$ of $G$ corresponding to $\lambda(\mathcal L)$, then $z_i$ is a constant for $i\in\mbox{sup}(\mathbf z)$ and $i$ being not the heart. Moreover, whenever $\mbox{sup}(\mathbf z)$ contains a vertex other than the heart, the signs of the heart and the vertices of leaves in $\mbox{sup}(\mathbf z)$ are opposite.
\end{Corollary}
However, in Section \ref{sec-odd}, we will show that $\mbox{sup}(\mathbf z)$ is a singleton which is the heart.

The next lemma is useful, which follows from a similar proof of \cite[Theorem 5]{q05}.

\begin{Lemma}\label{lem-5}
Let $k$ be even and $G=(V,E)$ be a $k$-uniform hypergraph. Let $\mathcal L$ be the Laplacian tensor of $G$. Then
\begin{eqnarray}\label{lp-c}
\lambda(\mathcal L)=\max\left\{\mathcal L\mathbf x^k:=\mathbf
x^T(\mathcal L\mathbf x^{k-1})\;|\;\sum_{i\in [n]}x_i^k=1,\;\mathbf
x\in\mathbb R ^n\right\}.
\end{eqnarray}
\end{Lemma}

The next lemma is an analogue of Corollary \ref{cor-4} for $k$ being even.
\begin{Lemma}\label{lem-2}
Let $k$ be even and $G=(V,E)$ be a hyperstar of size $d>0$. Then there is an H-eigenvector $\mathbf z\in\mathbb R^n$ of the Laplacian tensor $\mathcal L$ of $G$ satisfying that $z_i$ is a constant for $i\in\mbox{sup}(\mathbf z)$ and $i$ being not the heart.
\end{Lemma}

\noindent {\bf Proof.}
In the proof of Lemma \ref{lem-1}, $d>1$ is required only to guarantee $\lambda(\mathcal L)>1$. While, when $k$ is even, by Proposition \ref{prop-1}, $\lambda(\mathcal L)>1$ whenever $d>0$. Hence, there is an H-eigenvector $\mathbf x\in\mathbb R^n$ of the Laplacian tensor $\mathcal L$ of $G$ corresponding to $\lambda(\mathcal L)$ satisfying that $|x_i|$ is a constant for $i\in\mbox{sup}(\mathbf x)$ and $i$ being not the heart.

Suppose, without loss of generality, that $1$ is the heart. By Lemma \ref{lem-0}, without loss of generality, suppose that $\mbox{sup}(\mathbf x)=[n]$.
If $x_1>0$, then let $\mathbf y=-\mathbf x$, and otherwise let $\mathbf y=\mathbf x$.

Suppose that $y_i<0$ for some $i$ other than $y_1$. Then
\begin{eqnarray*}
0>(\lambda(\mathcal L)-1)y_i^{k-1}=-y_1\prod_{j\in le(i)\setminus\{1,i\}}y_j.
\end{eqnarray*}
Thus, a positive even number of vertices in $le(i)$ other than $1$ takes negative values.
Hence, all the values in this set
\begin{eqnarray*}
\left\{\prod_{j\in e_s\setminus\{1\}}y_j,\;e_s\in E\right\}
\end{eqnarray*}
are positive. Let $\mathbf z\in\mathbb R^n$ such that $z_1=y_1$ and $z_i=|y_i|$ for the others.
We have that if $y_i>0$, then
\begin{eqnarray*}
(\lambda(\mathcal L)-1)z_i^{k-1}=(\lambda(\mathcal L)-1)y_i^{k-1}=y_1\prod_{j\in le(i)\setminus\{1,i\}}y_j=z_1\prod_{j\in le(i)\setminus\{1,i\}}z_j;
\end{eqnarray*}
and if $y_i<0$, then
\begin{eqnarray*}
(\lambda(\mathcal L)-1)z_i^{k-1}=(\lambda(\mathcal L)-1)|y_i|^{k-1}=-y_1\prod_{j\in le(i)\setminus\{1,i\}}y_j=z_1\prod_{j\in le(i)\setminus\{1,i\}}z_j.
\end{eqnarray*}
Here, the second equality follows from the fact that $\prod_{j\in le(i)\setminus\{1,i\}}y_j<0$ in this situation.
Moreover,
\begin{eqnarray*}
(\lambda(\mathcal L)-d)z_1^{k-1}=(\lambda(\mathcal L)-d)y_1^{k-1}=\sum_{e_s\in E}\prod_{j\in e_s\setminus\{1\}}y_j=\sum_{e_s\in E}|\prod_{j\in e_s\setminus\{1\}}y_j|=\sum_{e_s\in E}\prod_{j\in e_s\setminus\{1\}}z_j.
\end{eqnarray*}
Consequently, $\mathbf z$ is the desired H-eigenvector. \ep

The next theorem gives the largest Laplacian H-eigenvalue of a hyperstar for $k$ being even.
\begin{Theorem}\label{thm-3}
Let $k$ be even and $G=(V,E)$ be a hyperstar of size $d>0$. Let $\mathcal L$ be the Laplacian tensor of $G$.
Then $\lambda(\mathcal L)$ is the unique real root of the equation $(1-\lambda)^{k-1}(\lambda-d)+d=0$ in the interval $(d,d+1)$.
\end{Theorem}

\noindent {\bf Proof.} By Lemma \ref{lem-2}, there is an H-eigenvector $\mathbf x\in\mathbb R^n$ of the Laplacian tensor $\mathcal L$ of $G$ satisfying that $x_i$ is a constant for $i\in\mbox{sup}(\mathbf x)$ and $i$ being not the heart.
By the proof for Lemma \ref{lem-0}, we have that $\lambda(\mathcal L)$ is the largest real root of the equation $(1-\lambda)^{k-1}(\lambda-w)+w=0$. Here $w$ is the size of the sub-hyperstar $G_{\mbox{sup}(\mathbf x)}$ of $G$.

Let $f(\lambda):=(1-\lambda)^{k-1}(\lambda-d)+d$. Then, $f^{\prime}(\lambda)=(1-\lambda)^{k-2}((k-1)(d-\lambda)+1-\lambda)$. Hence, $f$ is strictly decreasing in the interval $(d,+\infty)$. Moreover, $f(d+1)<0$. Consequently, $f$ has a unique real root in the interval $(d,d+1)$ which is the maximum. Thus, by Proposition \ref{prop-1}, we must have $\mbox{sup}(\mathbf x)=[n]$. The result follows. \ep

The next corollary is a direct consequence of Theorem \ref{thm-3}.
\begin{Corollary}\label{cor-1}
Let $G_1=(V_1,E_1)$ and $G_2=(V_2,E_2)$ be two hyperstars of size
$d_1$ and $d_2>0$, respectively. Let $\mathcal L_1$ and $\mathcal
L_2$ be the Laplacian tensors of $G_1$ and $G_2$ respectively. If
$d_1>d_2$, then $\lambda(\mathcal L_1)>\lambda(\mathcal L_2)$.
\end{Corollary}

When $k$ is even, the proofs of Lemmas \ref{lem-1} and \ref{lem-2}, and Theorem \ref{thm-3} actually imply the next corollary.
\begin{Corollary}\label{cor-3}
Let $k$ be even and $G=(V,E)$ be a hyperstar of size $d>0$. If $\mathbf x\in\mathbb R^n$ is an H-eigenvector of the Laplacian tensor $\mathcal L$ of $G$ corresponding to $\lambda(\mathcal L)$, then $\mbox{sup}(\mathbf x)=[n]$. Hence, there is an H-eigenvector $\mathbf z\in\mathbb R^n$ of the Laplacian tensor $\mathcal L$ of $G$ corresponding to $\lambda(\mathcal L)$ satisfying that $z_i$ is a constant for all the vertices other than the heart.
\end{Corollary}

\subsection{Even-Uniform Hypergraphs}

In this subsection, we present a tight lower bound for the largest
Laplacian H-eigenvalue and characterize the extreme hypergraphs when
$k$ is even.

The next theorem gives the lower bound, which is tight by Theorem
\ref{thm-3}.
\begin{Theorem}\label{thm-4}
Let $k$ be even and $G=(V,E)$ be a $k$-uniform hypergraph with the
maximum degree being $d>0$. Let $\mathcal L$ be the Laplacian tensor
of $G$. Then $\lambda(\mathcal L)$ is not smaller than the unique
real root of the equation $(1-\lambda)^{k-1}(\lambda-d)+d=0$ in the
interval $(d,d+1)$.
\end{Theorem}

\noindent {\bf Proof.} Suppose that $d_s=d$, the maximum degree.
Let $G'=(V',E')$ be a $k$-uniform hypergraph such that $E'=E_s$ and $V'$ consisting of the vertex $s$ and the vertices which share an edge with $s$. Let $\mathcal L'$ be the Laplacian tensor of $G'$. We claim that $\lambda(\mathcal L)\geq \lambda(\mathcal L')$.

Suppose that $|V'|=m\leq n$ and $\mathbf y\in\mathbb R^m$ is an H-eigenvector of $\mathcal L'$ corresponding to the H-eigenvalue $\lambda(\mathcal L')$ such that $\sum_{j\in [m]}y_j^k=1$. Suppose, without loss of generality, that $V'=[m]$, and the degree of vertex $j\in [m]$ in the hypergraph $G'$ is $d_j'$. Let $\mathbf x\in\mathbb R^n$ such that
\begin{eqnarray}\label{new-5}
x_i=y_i,\;\forall i\in [m],\;\mbox{and}\;x_i=0,\forall i>m.
\end{eqnarray}
Obviously, $\sum_{i\in [n]}x_i^k=\sum_{j\in [m]}y_j^k=1$. Moreover,
\begin{eqnarray}
\mathcal L\mathbf x^k&=&\sum_{i\in [n]}d_ix_i^k-k\sum_{e\in E}\prod_{j\in e}x_j\nonumber\\
&=&d_sx_s^k+\sum_{j\in [m]\setminus\{s\}}d'_jx_j^k-k\sum_{e\in E_s}\prod_{t\in e}x_t\nonumber\\
&&+\sum_{j\in [m]\setminus\{s\}}(d_j-d'_j)x_j^k+\sum_{j\in [n]\setminus [m]}d_jx_j^k-k\sum_{e\in E\setminus E_s}\prod_{t\in e}x_t\nonumber\\
&=&d_sx_s^k+\sum_{j\in [m]\setminus\{s\}}d'_jx_j^k-k\sum_{e\in E_s}\prod_{j\in e}x_j+\sum_{e\in E\setminus E_s}\left(\sum_{t\in e}x_t^k-k\prod_{w\in e}x_w\right)\nonumber\\
&=&\mathcal L'\mathbf y^k+\sum_{e\in E\setminus E_s}\left(\sum_{t\in e}x_t^k-k\prod_{w\in e}x_w\right)\nonumber\\
&\geq&\mathcal L'\mathbf y^k\label{ext-2}\\
&=&\lambda(\mathcal L').\nonumber
\end{eqnarray}
Here the inequality follows from the fact that $\sum_{t\in e}x_t^k-k\prod_{w\in e}|x_w|\geq 0$ by the arithmetic-geometric mean inequality. Thus, by the characterization \reff{lp-c} (Lemma \ref{lem-5}), we get the conclusion since $\lambda(\mathcal L)\geq \mathcal L\mathbf x^k$.

For the hypergraph $G'$, we define a new hypergraph by renumbering the vertices in the following way: fix the vertex $s$, and for every edge $e\in E_s$, number the rest $k-1$ vertices as $\{(e,2),\ldots,(e,k)\}$.
Let $\bar G=(\bar V,\bar E)$ be the $k$-uniform hypergraph with $\bar V:=\{s,(e,2),\ldots,(e,k),\;\forall e\in E_s\}$ and $\bar E:=\{\{s,(e,2),\ldots,(e,k)\}\;|\;e\in E_s\}$. It is easy to see that $\bar G$ is a hyperstar with size $d>0$ and the heart being $s$ (Definition \ref{def-bi-hm}). Let $\mathbf z\in\mathbb R^{kd-k+1}$ be an H-eigenvector of the Laplacian tensor $\bar \mathcal L$ of $\bar G$ corresponding to $\lambda(\bar\mathcal L)$. Suppose that $\sum_{t\in \bar V}z_t^k=1$. By Corollary \ref{cor-3}, we can choose a $\mathbf z$ such that $z_i$ is a constant other than $z_s$ which corresponds to the heart. Let $\mathbf y\in \mathbb R^m$ be defined as $y_i$ being the constant for all $i\in [m]\setminus \{s\}$ and $y_s=z_s$. Then, by a direct computation, we see that
\begin{eqnarray*}
\mathcal L'\mathbf y^k=\bar\mathcal L\mathbf z^k=\lambda(\bar\mathcal L).
\end{eqnarray*}
Moreover, $\sum_{j\in [m]}y_j^k\leq \sum_{t\in \bar V}z_t^k=1$. By \reff{lp-c} and the fact that $\lambda(\bar\mathcal L)>0$ (Theorem \ref{thm-3}), we see that
\begin{eqnarray}\label{ext}
\lambda(\mathcal L')\geq \lambda(\bar \mathcal L).
\end{eqnarray}
Consequently, $\lambda(\mathcal L)\geq \lambda(\bar\mathcal L)$. By
Theorem \ref{thm-3}, $\lambda(\bar \mathcal L)$ is the unique real
root of the equation $(1-\lambda)^{k-1}(\lambda-d)+d=0$ in the
interval $(d,d+1)$. Consequently, $\lambda(\mathcal L)$ is no
smaller than the unique real root of the equation
$(1-\lambda)^{k-1}(\lambda-d)+d=0$ in the interval $(d,d+1)$. \ep

By the proof of Theorem \ref{thm-4}, the next theorem follows
immediately.
\begin{Theorem}\label{cor-2}
Let $k$ be even, and $G=(V,E)$ and $G'=(V',E')$ be two $k$-uniform hypergraphs. Suppose that $\mathcal L$ and $\mathcal L'$ be the Laplacian tensors of $G$ and $G'$ respectively.
If $V\subseteq V'$ and $E\subseteq E'$, then $\lambda (\mathcal L)\leq\lambda(\mathcal L')$.
\end{Theorem}

The next lemma helps us to characterize the extreme hypergraphs with
respect to the lower bound of the largest Laplacian H-eigenvalue.
\begin{Lemma}\label{lem-4}
Let $k\geq 4$ be even and $G=(V,E)$ be a hyperstar of size $d>0$. Then there is an H-eigenvector $\mathbf z\in\mathbb R^n$ of the Laplacian tensor $\mathcal L$ of $G$ satisfying that exactly two vertices other than the heart in every edge takes negative values.
\end{Lemma}

\noindent {\bf Proof.} Suppose, without loss of generality, that $1$ is the heart. By Corollary \ref{cor-3}, there is an H-eigenvector $\mathbf x\in\mathbb R^n$ of $\mathcal L$ corresponding to $\lambda(\mathcal L)$ such that $x_i$ is a constant for the vertices other than the heart.
By Theorem \ref{thm-3}, we have that this constant is nonzero.
If $x_2<0$, then let $\mathbf y=-\mathbf x$, and otherwise let $\mathbf y=\mathbf x$. We have that $\mathbf y$ is an H-eigenvector of $\mathcal L$ corresponding to $\lambda(\mathcal L)$.

Let $\mathbf z\in\mathbb R^n$. We set $z_1=y_1$, and for every edge
$e\in E$ arbitrarily two chosen $i_{e,1},i_{e,2}\in e\setminus\{1\}$
we set $z_{i_{e,1}}=-y_{i_{e,1}}<0$, $z_{i_{e,2}}=-y_{i_2}<0$ and
$z_j=y_j>0$ for the others $j\in e\setminus\{1,i_{e,1},i_{e,2}\}$.
Then, by a direct computation, we can conclude that $\mathbf z$ is
an H-eigenvector of $\mathcal L$ corresponding to $\lambda(\mathcal
L)$. \ep

The next theorem is the main result of this subsection, which
characterizes the extreme hypergraphs with respect to the lower
bound of the largest Laplacian H-eigenvalue.
\begin{Theorem}\label{thm-5}
Let $k\geq 4$ be even and $G=(V,E)$ be a $k$-uniform connected hypergraph with the maximum degree being $d>0$. Let $\mathcal L$ be the Laplacian tensor of $G$.
Then $\lambda(\mathcal L)$ is equal to the unique real root of the equation $(1-\lambda)^{k-1}(\lambda-d)+d=0$ in the interval $(d,d+1)$ if and only if $G$ is a hyperstar.
\end{Theorem}

\noindent {\bf Proof.} By Theorem \ref{thm-3}, only necessity needs
a proof. In the following, suppose that $\lambda(\mathcal L)$ is
equal to the unique real root of the equation
$(1-\lambda)^{k-1}(\lambda-d)+d=0$ in the interval $(d,d+1)$.
Suppose that $d_s=d$ as before.

Define $G'$ and $\bar G$ as in Theorem \ref{thm-4}. Actually, let
$G'=(V',E')$ be the $k$-uniform hypergraph such that $E'=E_s$ and
$V'$ consisting of the vertex $s$ and the vertices which share an
edge with $s$. Let $\mathcal L'$ be the Laplacian tensor of $G'$.
Fix the vertex $s$, and for every edge $e\in E_s$, number the rest
$k-1$ vertices as $\{(e,2),\ldots,(e,k)\}$. Let $\bar G=(\bar V,\bar
E)$ be the $k$-uniform hypergraph such that $\bar
V:=\{s,(e,2),\ldots,(e,k),\;\forall e\in E_s\}$ and $\bar
E:=\{\{s,(e,2),\ldots,(e,k)\}\;|\;e\in E_s\}$.

With the same proof as in Theorem \ref{thm-4}, by Lemma \ref{lem-5},
we have that inequality in \reff{ext} is an equality if and only if
$|\bar V|=m$. Since otherwise $\sum_{j\in [m]}y_j^k<\sum_{t\in \bar
V}z_t^k=1$, which together with $\lambda(\bar \mathcal L)>0$ and
\reff{lp-c} implies that $\lambda(\mathcal L')>\lambda(\bar \mathcal
L)$. Hence, if $\lambda(\mathcal L)$ is equal to the unique real
root of the equation $(1-\lambda)^{k-1}(\lambda-d)+d=0$ in the
interval $(d,d+1)$, then $G'$ is a hyperstar. In this situation, the
inequality in \reff{ext-2} is an equality if and only if $G'=G$. The
sufficiency is clear.

For the necessity, suppose that $G'\neq G$. Then there is an edge
$\bar e\in E$
\begin{itemize}
\item [(i)] either containing both vertices in $[m]$ and vertices in $[n]\setminus [m]$, since $G$ is connected,
\item [(ii)] or containing only vertices in $[m]\setminus\{s\}$.
\end{itemize}
For the case (i), it is easy to get a contradiction since $\sum_{t\in e}x_t^k-k\prod_{w\in e}x_w=\sum_{t\in e\cap [m]}x_t^k>0$. Note that this situation happens if and only if $m<n$.
Then, in the following we assume that that $m=n$. For the case (ii), we must have that there are $q\geq 2$ edges $e_a\in E_s, a\in [q]$ in $G'$ such that $e_a\cap \bar e\neq\emptyset$ for all $a\in [q]$.
By Lemma \ref{lem-4}, let $\mathbf y\in\mathbb R^n$ be an H-eigenvector of the Laplacian tensor $\mathcal L'$ of $G'$ satisfying that exactly two vertices other than the heart in every edge takes negative values. Moreover, we can normalize $\mathbf y$ such that $\sum_{i\in [n]}y_i^k=1$. Since $m=n$, by \reff{new-5}, we have $\mathbf x=\mathbf y$. Consequently, by Lemma \ref{lem-5}, we have
\begin{eqnarray}
\lambda(\mathcal L)\geq
\mathcal L\mathbf x^k&=&\mathcal L'\mathbf x^k+\sum_{e\in E\setminus E_s}\left(\sum_{t\in e}x_t^k-k\prod_{w\in e}x_w\right)\nonumber\\
&=&\lambda(\mathcal L')+\sum_{e\in E\setminus E_s}\left(\sum_{t\in e}x_t^k-k\prod_{w\in e}x_w\right)\nonumber\\
&\geq &\lambda(\mathcal L')+\sum_{t\in \bar e}x_t^k-k\prod_{w\in \bar e}x_w.\nonumber
\end{eqnarray}
If $\prod_{w\in \bar e}x_w<0$, then we get a contradiction since $\lambda(\mathcal L')$ is equal to the unique real root of the equation $(1-\lambda)^{k-1}(\lambda-d)+d=0$ in the interval $(d,d+1)$.
In the following, we assume that $\prod_{w\in \bar e}x_w>0$. We have two cases:
\begin{itemize}
\item [(1)] $x_w>0$ or $x_w<0$ for all $w\in \bar e$,
\item [(2)] $x_b>0$ for some $b\in\bar e$ and $x_c<0$ for some $c\in\bar e$.
\end{itemize}
Note that $|e_a\cap \bar e|\leq k-2$ for all $a\in [q]$.
For an arbitrary but fixed $a\in [q]$, define $\{f_1,f_2\}:=\{f\in e_a\setminus\{s\}\;|\;x_f<0\}$.

(I). If $f_1, f_2 \in \bar e$, then we choose an $h\in e_a$ such that $h\neq s$, $h\notin \bar e$ and $x_h>0$. Since $k\geq 4$ is even, such an $h$ exists. It is a direct computation to see that $\mathbf z\in\mathbb R^n$ such that $z_{f_1}=-x_{f_1}>0$, $z_h=-x_h<0$, and $z_i=x_i$ for the others is still an H-eigenvector of $\mathcal L'$ corresponding to $\lambda(\mathcal L')$. More importantly, $\prod_{w\in \bar e}z_w<0$. Hence, replacing $\mathbf y$ by $\mathbf z$, we get a contradiction.

(II). If $f_1\in \bar e$ and $f_2\notin\bar e$, then either there is an $h\in \bar e\cap e_a$ such that $h\neq s$ and $x_h>0$, or there is an $h\in e_a$ such that $h\neq s$, $h\notin \bar e$ and $x_h>0$. Since $k\geq 4$ is even, such an $h$ exists. For the former case, set $\mathbf z\in\mathbb R^n$ such that $z_{h}=-x_{h}<0$, $z_{f_2}=-x_{f_2}>0$, and $z_i=x_i$ for the others; and for the latter case, set
$\mathbf z\in\mathbb R^n$ such that $z_{f_1}=-x_{f_1}>0$, $z_h=-x_h<0$, and $z_i=x_i$ for the others. Then, it is a direct computation to see that $\mathbf z$ is still an H-eigenvector of $\mathcal L'$ corresponding to $\lambda(\mathcal L')$. We also have that $\prod_{w\in \bar e}z_w<0$. Hence, replacing $\mathbf y$ by $\mathbf z$, we get a contradiction.

(III). The proof for the case $f_2\in \bar e$ and $f_1\notin\bar e$ is similar.

(IV). If $f_1,f_2\notin \bar e$, then there is some $b\in\bar e\cap e_a$ such that $x_b>0$, then similarly it is a direct computation to see that $\mathbf z\in\mathbb R^n$ such that $z_{b}=-x_{b}<0$, $z_{f_1}=-x_{f_1}>0$, and $z_i=x_i$ for the others is still an H-eigenvector of $\mathcal L'$ corresponding to $\lambda(\mathcal L')$. We also have that $\prod_{w\in \bar e}z_w<0$. Consequently, a contradiction can be derived.

Thus, $G=G'$ is a hyperstar.  \ep

Theorems \ref{thm-4} and \ref{thm-5} generalize the classical result
for graphs \cite{gm94,zl02}.

\subsection{Odd-Uniform Hypergraphs}\label{sec-odd}

In this subsection, we discuss odd-uniform hypergraphs. Note that there does not exist an analogue of Lemma \ref{lem-5} for $k$ being odd. Hence it is difficult to characterize the extreme hypergraphs for the lower bound of the largest H-eigenvalue of the Laplacian tensor.
\begin{Theorem}\label{thm-2}
Let $k$ be odd and $G=(V,E)$ be a hyperstar of size $d>0$. Let $\mathcal L$ be the Laplacian tensor of $G$.
Then $\lambda(\mathcal L)=d$.
\end{Theorem}

\noindent {\bf Proof.} The case for $d=1$ follows by direct computation, since in this case, for all $i\in [k]$
\begin{eqnarray*}
(\lambda(\mathcal L)-1)x_i^k=-\prod_{j\in [k]}x_j.
\end{eqnarray*}
If $\lambda(\mathcal L)>1$, then $x_i^k=x_j^k$ for all $i, j\in
[k]$. Since $k$ is odd and $x \not = 0$, we have $x_i=x_j \not = 0$
for all $i,j\in [k]$. This implies that $0<\lambda(\mathcal
L)-1=-1<0$, a contradiction.

In the following, we consider cases when $d>1$.
Suppose, without loss of generality, that $1$ is the heart. It is easy to see that the H-eigenvector $\mathbf x:=(1,0,\ldots,0)\in\mathbb R^n$ corresponds to the H-eigenvalue $d$. Suppose that $\mathbf x\in\mathbb R^n$ is an H-eigenvector of $\mathcal L$ corresponding to $\lambda(\mathcal L)$. In the following, we show that $\mbox{sup}(\mathbf x)=\{1\}$, which implies that $\lambda(\mathcal L)=d$.

Suppose on the contrary that $\mbox{sup}(\mathbf x)\neq\{1\}$. By Lemma \ref{lem-0} and Corollary \ref{cor-4}, without loss of generality, we assume that $\mbox{sup}(\mathbf x)=[n]$ and $\mathbf x$ is of the following form
\begin{eqnarray*}
\alpha:=x_1>0, \;\mbox{and}\; x_2=\cdots=x_n=-1.
\end{eqnarray*}
Then, we see that
\begin{eqnarray*}
(\mathcal L\mathbf x^{k-1})_1=d \alpha^{k-1}-d=\lambda(\mathcal L)\alpha^{k-1},
\end{eqnarray*}
and for $i\in\{2,\ldots,n\}$
\begin{eqnarray*}
(\mathcal L\mathbf x^{k-1})_i=1+\alpha=\lambda(\mathcal L).
\end{eqnarray*}
Consequently,
\begin{eqnarray*}
(d-\lambda(\mathcal L))(\lambda(\mathcal L)-1)^{k-1}=d.
\end{eqnarray*}
Hence, we must have $\lambda(\mathcal L)<d$. This is a contradiction. Hence, $\lambda(\mathcal L)=d$. \ep

When $k$ is odd, Theorem \ref{thm-2}, together with Lemma
\ref{lem-00}, implies that the maximum degree is a tight lower bound
for the largest Laplacian H-eigenvalue.

We now give a lower bound for the largest Laplacian H-eigenvalue of
a $3$-uniform complete hypergraph.
\begin{Proposition}\label{prop-2}
Let $G=(V,E)$ be a $3$-uniform complete hypergraph. Let $\mathcal L$
be the Laplacian tensor of $G$ and $n=2m$ for some positive integer
$m$. Then $\lambda(\mathcal L)\geq {n-1\choose 2}+m-1$, which is
strictly larger than $d = {n-1\choose 2}$, the maximum degree of
$G$.
\end{Proposition}

\noindent {\bf Proof.} Let $\mathbf x\in\mathbb R^n$ be defined as
$x_1=\cdots=x_m=1$ and $x_{m+1}=\cdots=x_{2m}=-1$. We have that
\begin{eqnarray*}
\left(\mathcal L\mathbf x^2\right)_1&=&{n-1\choose 2}x_1^2 -\sum_{1<i<j\in [n]}x_ix_j\\
&=&{n-1\choose 2}-\sum_{1<i<j\in [m]}x_ix_j-\sum_{m+1\leq i<j\in [2m]}x_ix_j-\sum_{1<i\in [m],\;m+1\leq j\leq 2m}x_ix_j\\
&=&{n-1\choose 2}-\sum_{1<i<j\in [m]}x_ix_j-\sum_{m+1\leq i<j\in [2m]}x_ix_j+\sum_{1<i\in [m],\;m+1\leq j\leq 2m}\left|x_ix_j\right|\\
&=&{n-1\choose 2}-{m-1\choose 2}-{m\choose 2}+(m-1)m\\
&=&\left[{n-1\choose 2}+m-1\right]x_1^2.
\end{eqnarray*}
Thus, for any $p = 2, \cdots, m$, we have that
$$\left(\mathcal L\mathbf x^2\right)_p = \left(\mathcal L\mathbf
x^2\right)_1 = \left[{n-1\choose 2}+m-1\right]x_p^2.$$ Similarly,
for any $p\in \{m+1,\ldots, 2m\}$, we have that
\begin{eqnarray*}
\left(\mathcal L\mathbf x^2\right)_p&=&\left(\mathcal L\mathbf
x^2\right)_n\\
&=&{n-1\choose 2}x_n^2 -\sum_{1\leq i<j\in [n-1]}x_ix_j\\
&=&{n-1\choose 2}-\sum_{1\leq i<j\in [m]}x_ix_j-\sum_{m+1\leq i<j\in [2m-1]}x_ix_j-\sum_{1\leq i\in [m],\;m+1\leq j\leq 2m-1}x_ix_j\\
&=&{n-1\choose 2}-\sum_{1\leq i<j\in [m]}x_ix_j-\sum_{m+1\leq i<j\in [2m-1]}x_ix_j+\sum_{1\leq i\in [m],\;m+1\leq j\leq 2m-1}\left|x_ix_j\right|\\
&=&{n-1\choose 2}-{m\choose 2}-{m-1\choose 2}+m(m-1)\\
&=&\left[{n-1\choose 2}+m-1\right]x_p^2.
\end{eqnarray*}
Thus, $\mathbf x$ is an H-eigenvector of $\mathcal L$ corresponding
to the H-eigenvalue ${n-1\choose 2}+m-1$. \ep

We have the following conjecture.
\begin{Conjecture}\label{con-1}
Let $k\geq 3$ be odd and $G=(V,E)$ be a $k$-uniform connected hypergraph with the maximum degree being $d>0$. Let $\mathcal L$ be the Laplacian tensor of $G$.
Then $\lambda(\mathcal L)$ is equal to $d$ if and only if $G$ is a hyperstar.
\end{Conjecture}

\section{The Largest Signless Laplacian H-eigenvalue}
\setcounter{Theorem}{0} \setcounter{Proposition}{0}
\setcounter{Corollary}{0} \setcounter{Lemma}{0}
\setcounter{Definition}{0} \setcounter{Remark}{0}
\setcounter{Conjecture}{0}  \setcounter{Example}{0} \hspace{4mm} In
this section, we discuss the largest signless Laplacian H-eigenvalue
of a $k$-uniform hypergraph. Since the signless Laplacian tensor
$\mathcal Q$ is nonnegative, the situation is much clearer than the
largest Laplacian H-eigenvalue.

The next proposition gives bounds on $\lambda(\mathcal Q)$.
\begin{Proposition}\label{thm-41}
Let $G=(V,E)$ be a $k$-uniform hypergraph with maximum degree being $d>0$, and $\mathcal A$ and $\mathcal{Q}$ be the adjacency tensor and the signless Laplacian tensor of $G$ respectively. Then
\begin{eqnarray*}
\max\left\{d,\frac{2\sum_{i\in [n]}d_i}{n}\right\}\leq
\lambda(\mathcal{Q})\leq \lambda(\mathcal{A})+d.
\end{eqnarray*}
\end{Proposition}

\noindent {\bf Proof.} The first inequality follows from \cite[Corollary 12]{q12a}. For the second, by \cite[Theorem 11]{q12a}, we have that
\begin{eqnarray*}
\lambda(\mathcal{Q})&=&\max_{\sum_{i\in [n]}x_i^k=1,\;\mathbf x\in\mathbb R_+^n} \mathcal{Q}\mathbf{x}^k=\max_{\sum_{i\in [n]}x_i^k=1,\;\mathbf x\in\mathbb R_+^n} (\mathcal{A}+\mathcal{D})\mathbf{x}^k\nonumber\\
&\leq&\max_{\sum_{i\in [n]}x_i^k=1,\;\mathbf x\in\mathbb R_+^n} \mathcal{A}\mathbf{x}^k+\max_{\sum_{i\in [n]}x_i^k=1,\;\mathbf x\in\mathbb R_+^n} \mathcal{D}\mathbf{x}^k=\lambda(\mathcal{A})+d.
\end{eqnarray*}
Consequently, the second inequality follows. \ep

\begin{Lemma}\label{lem-6}
Let $G =(V,E)$ be a $k$-uniform regular connected hypergraph with degree $d>0$, and $\mathcal Q$ be its signless Laplacian tensor. Then, $\lambda(\mathcal Q)=2d$.
\end{Lemma}

\noindent {\bf Proof.} Note that the vector of all ones is an H-eigenvector of $\mathcal Q$ corresponding to the H-eigenvalue $2d$. Since $\mathcal Q$ is weakly irreducible (\cite[Lemma 3.1]{pz12}), the result follows from \cite[Lemmas 2.2 and 2.3]{hq12a}. \ep

The next proposition gives a tight upper bound of the largest signless Laplacian H-eigenvalues and characterizes the extreme hypergraphs.

\begin{Proposition}\label{thm-42}
Let $G =(V,E)$ be a $k$-uniform hypergraph and $G'$ be a sub-hypergraph of $G$. Let $\mathcal{Q}$ and $\mathcal{Q}'$ be the signless Laplacian tensor of $G$ and $G'$, respectively. Then,
\begin{eqnarray} \lambda(\mathcal{Q}') \leq \lambda(\mathcal{Q}).\nonumber\end{eqnarray}
Furthermore, if $G'$ and $G$ are both connected, then $\lambda(\mathcal{Q}')=\lambda(\mathcal{Q})$ if and only if
$G'=G$. Consequently,
\begin{eqnarray}
\lambda(\mathcal{Q})&\leq & 2{n-1 \choose k-1},\nonumber
\end{eqnarray}
and equality holds if and only if $G$ is a $k$-uniform complete hypergraph.
\end{Proposition}

\noindent {\bf Proof.} The first conclusion follows from \cite[Theorem 3.19]{yy11}.
The remaining follows from \cite[Theorem 3.20]{yy11}, \cite[Theorem 4]{q12b} and \cite[Lemma 3.1]{pz12} (see also \cite[Lemmas 2.2 and 2.3]{hq13}) which imply that there is a unique positive H-eigenvector of $\mathcal Q$ and the corresponding H-eigenvalue must be $\lambda(\mathcal Q)$ whenever $G$ is connected, and the fact that
the vector of all ones is an H-eigenvector of $\mathcal Q$ corresponding to the H-eigenvalue $2{n-1 \choose k-1}$ when $G$ is a complete hypergraph (Lemma \ref{lem-6}).  \ep

When $k=2$ (i.e., the usual graph), Propositions \ref{thm-41} and \ref{thm-42} reduce to the classic results in graph theory \cite{crs07}.

The next theorem gives a tight lower bound for $\lambda(\mathcal Q)$ and characterizes the extreme hypergraphs.
\begin{Theorem}\label{thm-44}
Let $G=(V,E)$ be a $k$-uniform connected hypergraph with the maximum degree being $d>0$ and $\mathcal{Q}$ be the signless Laplician tensor of $G$. Then
\begin{eqnarray*}
\lambda(\mathcal{Q})&\geq &d+ d
\left(\frac{1}{\alpha_*}\right)^{k-1},
\end{eqnarray*}
where $\alpha_*\in (d-1,d]$ is the largest real root of
$\alpha^k+(1-d)\alpha^{k-1}-d=0$, with equality holding if and only if $G$ is a hyperstar.
\end{Theorem}

\noindent {\bf Proof.}  Suppose that $d_s=d$.
Let $G'$ be the hypergraph $G_{S}$ with $S$ being the vertices in the set $E_s$. As in the proof of Proposition \ref{thm-4}, for the hypergraph $G'$, we define a new hypergraph by renumbering the vertices in the following way: fix the vertex $s$, and for every edge $e\in E_s$, number the rest $k-1$ vertices as $\{(e,2),\ldots,(e,k)\}$.
Let $\bar G=(\bar V,\bar E)$ be the $k$-uniform hypergraph such that $\bar V:=\{s,(e,2),\ldots,(e,k),\;\forall e\in E_s\}$ and $\bar E:=\{\{s,(e,2),\ldots,(e,k)\}\;|\;e\in E_s\}$. It is easy to see that $\bar G$ is a hyperstar with order $d>0$ and the heart being $s$. Let $\mathbf z\in\mathbb R^{kd-k+1}$ be a vector such that $z_s=\alpha>0$ and $z_j=1$ for all $j\in\bar V\setminus\{s\}$. By a similar proof of Proposition \ref{prop-1}, we see that $\mathbf z$ is an H-eigenvector of the signless Laplacian tensor $\bar \mathcal Q$ of $\bar G$ if and only if $\alpha$ is a real root of the following equation
\begin{eqnarray}\label{sl}
\alpha^k+(1-d)\alpha^{k-1}-d=0.
\end{eqnarray}
In this situation, the H-eigenvalue is $\lambda=1+\alpha$.

By \cite[Theorem 11]{q12a} and \cite[Lemmas 2.2 and 2.3]{hq13}, or by a similar proof of Proposition \ref{prop-1}, we can show that $\lambda(\mathcal Q)\geq \lambda(\bar\mathcal Q)$ with equality holding if and only if $G=\bar G$.
Moreover, let $\alpha_*$ be the largest real root of the equation \reff{sl}, by \reff{sl} we have
\begin{eqnarray*}
\lambda_*=1+\alpha_*=d+d\left(\frac{1}{\alpha_*}\right)^{k-1}.
\end{eqnarray*}
With a similar proof as Theorem \ref{thm-3}, we can show that the equation in \reff{sl} has a unique real in the interval $(d-1, d]$ which is the maximum.
Since $\bar G$ is connected, by \cite[Lemmas 2.2 and 2.3]{hq13} and \cite[Lemma 3.1]{pz12}, we have that $\lambda(\bar\mathcal Q)=1+\alpha_*$. Consequently, the results follow. \ep

When $G$ is a 2-uniform hypergraph, we know that $\alpha_*=d$, hence
Theorem \ref{thm-44} reduces to $\lambda(\mathcal{Q})\geq d+1$
\cite{crs07}.

\section{The Relation between The Largest Laplacian and Signless Laplacian H-Eigenvalues}
\setcounter{Theorem}{0} \setcounter{Proposition}{0}
\setcounter{Corollary}{0} \setcounter{Lemma}{0}
\setcounter{Definition}{0} \setcounter{Remark}{0}
\setcounter{Conjecture}{0}  \setcounter{Example}{0} \hspace{4mm} In
this section, we discuss the relationship between the largest
Laplacian H-eigenvalue and the largest signless Laplacian
H-eigenvalue.

The following theorem characterizes this relationship.  This theorem
generalizes the classical result in spectral graph theory
\cite{zl02,z07}.

\begin{Theorem}\label{prop-3}
Let $G=(V,E)$ be a $k$-uniform hypergraph. Let $\mathcal L,\mathcal
Q$ be the Laplacian and signless Laplacian tensors of $G$
respectively. Then
\begin{eqnarray*}
\lambda(\mathcal L)\leq\lambda(\mathcal Q).
\end{eqnarray*}
If furthermore $G$ is connected and $k$ is even, then
\begin{eqnarray*}
\lambda(\mathcal L)=\lambda(\mathcal Q)
\end{eqnarray*}
if and only if $G$ is odd-bipartite.
\end{Theorem}

\noindent {\bf Proof.}  The first conclusion follows from Definition
\ref{def-00} and \cite[Proposition 14]{q12a}.

We now prove the second conclusion.   We first prove the
sufficiency. We assume that $G$ is odd-bipartite. Suppose that
$\mathbf x\in\mathbb R^n$ is a nonnegative H-eigenvector of
$\mathcal Q$ corresponding to $\lambda(\mathcal Q)$. Then,
\cite[Lemma 2.2]{hq12a} implies that $\mathbf x$ is a positive
vector, i.e., all its entries are positive.  Suppose that $V=V_1\cup
V_2$ is an odd-bipartition of $V$ such that $V_1,V_2\neq \emptyset$
and every edge in $E$ intersects $V_1$ with exactly an odd number of
vertices. Let $\mathbf y\in\mathbb R^n$ be defined such that
$y_i=x_i$ whenever $i\in V_1$ and $y_i=-x_i$ for the others. Then,
for $i\in V_1$, we have
\begin{eqnarray*}
\left[(\mathcal D-\mathcal A)\mathbf y^{k-1}\right]_i&=&d_iy_i^{k-1}-\sum_{e\in E_i}\prod_{j\in e\setminus\{i\}}y_j\\
&=&d_ix_i^{k-1}+\sum_{e\in E_i}\prod_{j\in e\setminus\{i\}}x_j\\
&=&\left[(\mathcal D+\mathcal A)\mathbf x^{k-1}\right]_i\\
&=&\lambda(\mathcal Q) x_i^{k-1}\\
&=&\lambda(\mathcal Q) y_i^{k-1}.
\end{eqnarray*}
Here the second equality follows from the fact that exactly an odd number of vertices in $e$ takes negative values for every $e\in E_i$.
Similarly, we have for $i\in V_2$,
\begin{eqnarray*}
\left[(\mathcal D-\mathcal A)\mathbf y^{k-1}\right]_i&=&d_iy_i^{k-1}-\sum_{e\in E_i}\prod_{j\in e\setminus\{i\}}y_j\\
&=&-d_ix_i^{k-1}-\sum_{e\in E_i}\prod_{j\in e\setminus\{i\}}x_j\\
&=&-\left[(\mathcal D+\mathcal A)\mathbf x^{k-1}\right]_i\\
&=&-\lambda(\mathcal Q) x_i^{k-1}\\
&=&\lambda(\mathcal Q) y_i^{k-1}.
\end{eqnarray*}
Here the second equality follows from the fact that exactly an even number of vertices in $e\setminus\{i\}$ takes negative values for every $e\in E_i$, and the last from the fact that $y_i=-x_i$. Thus, $\lambda(\mathcal Q)$ is an H-eigenvalue of $\mathcal L$. This, together with the first conclusion, implies that $\lambda(\mathcal L)=\lambda(\mathcal Q)$.

In the following, we prove the necessity of the second conclusion.
We assume that $\lambda(\mathcal L)=\lambda(\mathcal Q)$. Let
$\mathbf x\in\mathbb R^n$ be an H-eigenvector of $\mathcal L$
corresponding to the H-eigenvalue $\lambda(\mathcal L)$ such that
$\sum_{i\in [n]}x_i^k=1$. Then,
\begin{eqnarray*}
\left[(\mathcal D-\mathcal A)\mathbf x^{k-1}\right]_i=\lambda(\mathcal L)x_i^{k-1},\;\forall i\in [n].
\end{eqnarray*}
Let $\mathbf y\in\mathbb R^n$ be defined such that $y_i=|x_i|$ for all $i\in [n]$. By \reff{lp-c}, we see that
\begin{eqnarray}\label{new-1}
\lambda(\mathcal L)&=&\sum_{i\in [n]}x_i \left[(\mathcal D-\mathcal A)\mathbf x^{k-1}\right]_i=\sum_{i\in [n]}|x_i|\left| \left[(\mathcal D-\mathcal A)\mathbf x^{k-1}\right]_i\right|\nonumber\\
&\leq& \sum_{i\in [n]}y_i\left[(\mathcal D+\mathcal A)\mathbf y^{k-1}\right]_i\leq\lambda(\mathcal Q).
\end{eqnarray}
Thus, all the inequalities in \reff{new-1} should be equalities. By \cite[Lemma 2.2]{q12a} and \cite[Theorem 2.1(iii)]{hq13}, we have that $\mathbf y$ is an H-eigenvector of $\mathcal Q$ corresponding to the H-eigenvalue $\lambda(\mathcal Q)$, and it is a positive vector.
Let $V_1:=\{i\in [n]\;|\;x_i>0\}$ and $V_2:=\{i\in [n]\;|\;x_i<0\}$. Then, $V_1\cup V_2=[n]$, since $\mathbf y$ is positive. Since $G$ is connected and nontrivial, we must have that $V_2\neq \emptyset$. Otherwise $\left|\left[(\mathcal D-\mathcal A)\mathbf x^{k-1}\right]_i\right|<\left[(\mathcal D+\mathcal A)\mathbf y^{k-1}\right]_i$, since $(\mathcal A\mathbf x^{k-1})_i>0$ in this situation. We also have that $V_1\neq\emptyset$, since otherwise $\left|\left[(\mathcal D-\mathcal A)\mathbf x^{k-1}\right]_i\right|=|-d_iy_i^{k-1}+(\mathcal A\mathbf y^{k-1})_i|<\left[(\mathcal D+\mathcal A)\mathbf y^{k-1}\right]_i$.

Moreover, since the first inequality in \reff{new-1} must be an equality, we must get that for all $i\in V_1$,
\begin{eqnarray*}
\lambda(\mathcal Q)y_i^{k-1}&=&\left[(\mathcal D+\mathcal A)\mathbf y^{k-1}\right]_i=\left[(\mathcal D-\mathcal A)\mathbf x^{k-1}\right]_i
\end{eqnarray*}
We have that
\begin{eqnarray*}
\left[(\mathcal D+\mathcal A)\mathbf y^{k-1}\right]_i=d_iy_i^{k-1}+\sum_{e\in E_i}\prod_{j\in e\setminus\{i\}}y_j,\;\mbox{and}\;
\left[(\mathcal D-\mathcal A)\mathbf x^{k-1}\right]_i=d_ix_i^{k-1}-\sum_{e\in E_i}\prod_{j\in e\setminus\{i\}}x_j.
\end{eqnarray*}
Hence, for every $e\in E_i$ with $i\in V_1$, we must have that exactly $|e\cap V_2|$ is an odd number. Similarly, we can show that for every $e\in E_i$ with $i\in V_2$, we must have that exactly $|e\cap V_1|$ is an odd number. Consequently, $G$ is odd-bipartite by Definition \ref{def-bi-odd}.\ep

In the following, we give an application of Theorem \ref{prop-3}.

\begin{Definition}\label{def-hc}
Let $G=(V,E)$ be a $k$-uniform nontrivial
hypergraph. If there is a disjoint partition of the vertex set $V$ as $V=V_1\cup\cdots\cup V_s$ such that $|V_1|=\cdots=|V_s|=k$, and
\begin{itemize}
\item [(i)] $E=\{V_i\;|\;i\in [s]\}$,
\item [(ii)] $|V_1\cap V_2|=\cdots=|V_{s-1}\cap V_s|=|V_s\cap V_1|=1$, and $V_i\cap V_j=\emptyset$ for the other cases,
\item [(iii)] the intersections $V_1\cap V_2$, \ldots, $V_s\cap V_1$ are mutually different.
\end{itemize}
then $G$ is called a {\em hypercycle}. $s$ is the {\em size} of the hypercycle.
\end{Definition}
It is easy to see that a $k$-uniform hypercycle of size $s>0$ has $n=s(k-1)$ vertices, and is connected.
Figure 3 (i) is an example of a $4$-uniform hypercycle of size $3$.

\begin{figure}[htbp]
\centering
\includegraphics[width=2.0in]{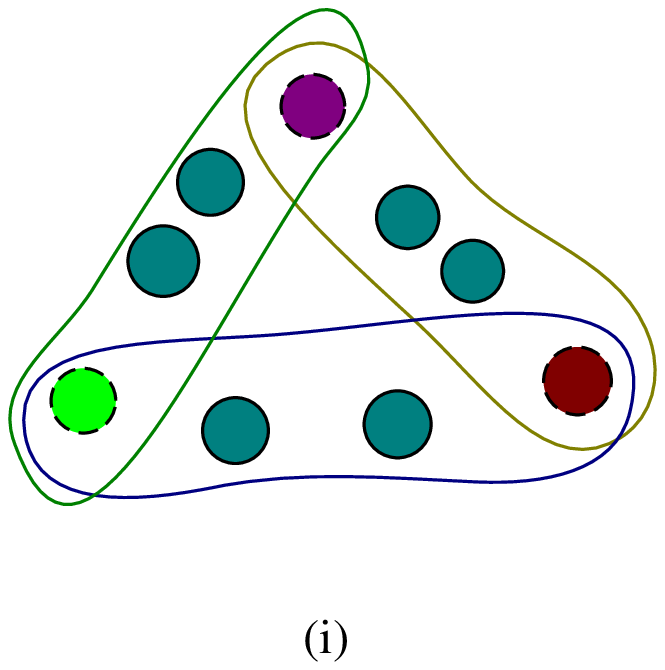}
\includegraphics[width=2.0in]{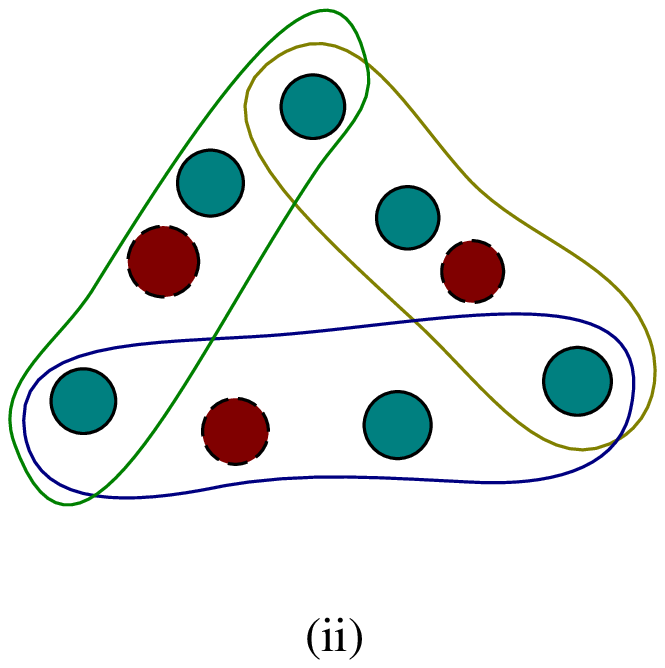}
\caption{(i) is an example of a $4$-uniform hypercycle of size $3$. The intersections are in dashed margins. (ii) is an illustration of an odd-bipartition of the $4$-uniform hypercycle. The partition is clear from the different colors of the disks (also the dashed margins from the solid ones).}
\end{figure}

The next lemma says that the largest signless Laplacian H-eigenvalue of a hypercycle is easy to characterize.
\begin{Lemma}\label{lem-7}
Let $G=(V,E)$ be a $k$-uniform
hypercycle of size $s>0$ and $\mathcal Q$ be its signless Laplacian tensor. Then, $\lambda(\mathcal Q)=2+2\beta^{k-2}$ with $\beta$ being the unique positive solution of the equation $2\beta^k+\beta^2-1=0$ which is in the interval $(\frac{1}{2}, 1)$.
\end{Lemma}

\noindent {\bf Proof.} By \cite[Theorem 3.20]{yy11}, \cite[Theorem 4]{q12b} and \cite[Lemma 3.1]{pz12} (see also \cite[Lemmas 2.2 and 2.3]{hq13}), if we can find a positive H-eigenvector $\mathbf x\in\mathbb R^n$ of $\mathcal Q$ corresponding to an H-eigenvalue $\mu$, then $\mu=\lambda(\mathcal Q)$.

Let $x_i=\alpha$ whenever $i$ is an intersection of the edges of $G$ and $x_i=\beta$ for the others. Without loss of generality, we assume that $\alpha =1$. Then, for an intersection vertex $i$, we have that $d_i=2$ and
\begin{eqnarray*}
(\mathcal Q\mathbf x^{k-1})_i=2\alpha^{k-1}+2\alpha\beta^{k-2}=2+2\beta^{k-2};
\end{eqnarray*}
and
for the other vertices $j$, we have that $d_j=1$ and
\begin{eqnarray*}
(\mathcal Q\mathbf x^{k-1})_j=\beta^{k-1}+\alpha^2\beta^{k-3}=\beta^{k-1}+\beta^{k-3}.
\end{eqnarray*}
If there are some $\mu>0$ and $\beta>0$ such that
\begin{eqnarray}\label{new-6}
2+2\beta^{k-2}=\mu,\;\mbox{and}\;\beta^{k-1}+\beta^{k-3}=\mu\beta^{k-1},
\end{eqnarray}
then $\mu=\lambda(\mathcal Q)$ by the discussion at the beginning of this proof.
We assume that \reff{new-6} has a required solution pair. Then,
\begin{eqnarray*}
2\beta^{2k-3}+\beta^{k-1}-\beta^{k-3}=0,\;i.e.,\; 2\beta^k+\beta^2-1=0.
\end{eqnarray*}
Let $g(\beta):=2\beta^k+\beta^2-1$. Then $g(1)>0$ and
\begin{eqnarray*}
g(\frac{1}{2})=\frac{1}{2^{k-1}}+\frac{1}{4}-1<0.
\end{eqnarray*}
Thus, \reff{new-6} does have a solution pair with $\beta\in (\frac{1}{2}, 1)$ and $\mu=2+2\beta^{k-2}$. Since $\mathcal Q$ has a unique positive H-eigenvector (\cite[Lemmas 2.2 and 2.3]{hq13}), the equation $2\beta^k+\beta^2-1=0$ has a unique positive solution which is in the interval $(\frac{1}{2}, 1)$. Hence, the result follows. \ep

By Theorem \ref{prop-3} and Lemma \ref{lem-7}, we can get the following corollary, which characterizes the largest Laplacian H-eigenvalue of a hypercycle when $k$ is even.
\begin{Corollary}\label{cor-5}
Let $k$ be even and $G=(V,E)$ be a $k$-uniform
hypercycle of size $s>0$. Let $\mathcal L$ be its Laplacian tensor. Then, $\lambda(\mathcal L)=2+2\beta^{k-2}$ with $\beta$ being the unique positive solution of the equation $2\beta^k+\beta^2-1=0$ which is in the interval $(\frac{1}{2}, 1)$.
\end{Corollary}

\noindent {\bf Proof.} By Theorem \ref{prop-3} and Lemma \ref{lem-7}, it suffices to show that when $k$ is even, a $k$-uniform
hypercycle is odd-bipartite.

Let $V=V_1\cup\cdots\cup V_s$ such that $|V_1|=\cdots=|V_s|=k$ be the partition of the vertices satisfying the hypotheses in Definition \ref{def-hc}.
Denote $V_s\cap V_1$ as $i_1$, $V_1\cap V_2$ as $i_2$, \ldots, $V_{s-1}\cap V_s$ as $i_s$. For every $j\in [s]$, choose a vertex $r_j\in V_j$ such that $r_j\notin\{i_1,\ldots,i_s\}$. Let $S_1:=\{r_j\;|\;j\in [s]\}$ and $S_2=V\setminus S_1$.
Then it is easy to see that $S_1\cup S_2=V$ is an odd-bipartition of $G$ (Definition \ref{def-bi-odd}). An illustration of such a partition is shown in Figure 3 (ii).

Thus, the result follows. \ep

The next proposition says that when $k$ is odd, the two H-eigenvalues cannot equal for a connected nontrivial hypergraph.

\begin{Proposition}\label{prop-4}
Let $k$ be odd and $G=(V,E)$ be a $k$-uniform connected nontrivial
hypergraph. Let $\mathcal L,\mathcal Q$ be the Laplacian and signless Laplacian tensors of $G$ respectively.
Then
\begin{eqnarray*}
\lambda(\mathcal L)<\lambda(\mathcal Q).
\end{eqnarray*}
\end{Proposition}

\noindent {\bf Proof.} Suppose that $\mathbf x\in\mathbb R^n$ is an H-eigenvector of $\mathcal L$ corresponding to $\lambda(\mathcal L)$ such that $\sum_{i\in [n]}|x_i|^k=1$. Then, we have that
\begin{eqnarray*}
\lambda(\mathcal L)x_i^{k-1}=(\mathcal L\mathbf x^{k-1})_i=\left[(\mathcal D-\mathcal A)\mathbf x^{k-1}\right]_i,\;\forall i\in [n].
\end{eqnarray*}
Hence,
\begin{eqnarray}\label{new-2}
\lambda(\mathcal L)&=&\sum_{i\in [n]}|x_i|\left|(\mathcal L\mathbf x^{k-1})_i\right|=\sum_{i\in [n]}|x_i|\left|\left[(\mathcal D-\mathcal A)\mathbf x^{k-1}\right]_i\right|\nonumber\\
&\leq& \sum_{i\in [n]}|x_i|\left[(\mathcal D+\mathcal A)|\mathbf x|^{k-1}\right]_i\leq \lambda(\mathcal Q).
\end{eqnarray}
If $\mbox{sup}(\mathbf x)\neq [n]$, then $\lambda(\mathcal L)<\lambda(\mathcal Q)$ by \cite[Lemma 2.2]{hq13}. Hence, in the following we assume that $\mbox{sup}(\mathbf x)=[n]$. We prove the conclusion by contradiction. Suppose that $\lambda(\mathcal L)=\lambda(\mathcal Q)$. Then all the inequalities in \reff{new-2} should be equalities. By \cite[Theorem 11]{q12a},  $\mathbf y:=|\mathbf x|$ is an H-eigenvector of $\mathcal Q$ corresponding to the H-eigenvalue $\lambda(\mathcal Q)$, and it is a positive vector. Similar to the proof of Proposition \ref{prop-3}, we can get a bipartition of $V$ as $V=V_1\cup V_2$ with $V_1,V_2\neq\emptyset$. Moreover, for all $i\in V$,
\begin{eqnarray*}
\lambda(\mathcal Q)y_i^{k-1}&=&\left[(\mathcal D+\mathcal A)\mathbf y^{k-1}\right]_i=\left|\left[(\mathcal D-\mathcal A)\mathbf x^{k-1}\right]_i\right|.
\end{eqnarray*}
Suppose, without loss of generality, that $x_1>0$. Then, we have that $|e\cap V_2|<k-1$ is an odd number for every $e\in E_1$. Since $G$ is connected and nontrivial, we have that $E_1\neq \emptyset$. Suppose that $2\in\bar e\cap V_2$ with $\bar e\in E_1$. We have $x_2<0$ and
\begin{eqnarray*}
\left|\left[(\mathcal D-\mathcal A)\mathbf x^{k-1}\right]_2\right|&=&\left|d_2x_2^{k-1}-\sum_{e\in E_2}\prod_{w\in e}x_w\right|\\
&=&\left|d_2x_2^{k-1}-\sum_{e\in E_2\setminus\{\bar e\}}\prod_{w\in e}x_w-\prod_{w\in\bar e}x_w\right|\\
&=&\left|d_2|x_2|^{k-1}-\sum_{e\in E_2\setminus\{\bar e\}}\prod_{w\in e}x_w-\prod_{w\in\bar e}|x_w|\right|\\
&\leq&\left|\left|d_2|x_2|^{k-1}+\sum_{e\in E_2\setminus\{\bar e\}}\prod_{w\in e}|x_w|\right|-\prod_{w\in\bar e}|x_w|\right|\\
&<&\left|d_2|x_2|^{k-1}+\sum_{e\in E_2}\prod_{w\in e}|x_w|\right|\\
&=&\left[(\mathcal D-\mathcal A)\mathbf y^{k-1}\right]_2.
\end{eqnarray*}
Thus, we get a contradiction. Consequently, $\lambda(\mathcal L)<\lambda(\mathcal Q)$. \ep

Combining Theorem \ref{prop-3} and Proposition \ref{prop-4}, we have
the following theorem.

\begin{Theorem}
Let $G=(V,E)$ be a $k$-uniform hypergraph. Let $\mathcal L,\mathcal
Q$ be the Laplacian and signless Laplacian tensors of $G$
respectively. Then
\begin{eqnarray*}
\lambda(\mathcal L)\leq\lambda(\mathcal Q).
\end{eqnarray*}
If furthermore $G$ is connected, then
\begin{eqnarray*}
\lambda(\mathcal L)=\lambda(\mathcal Q)
\end{eqnarray*}
if and only if $k$ is even and $G$ is odd-bipartite.
\end{Theorem}

\section{Final Remarks}
\setcounter{Theorem}{0} \setcounter{Proposition}{0}
\setcounter{Corollary}{0} \setcounter{Lemma}{0}
\setcounter{Definition}{0} \setcounter{Remark}{0}
\setcounter{Conjecture}{0}  \setcounter{Example}{0} \hspace{4mm} In
this paper, the largest Laplacian and signless Laplacian
H-eigenvalues of a uniform hypergraph are discussed. The largest
signless Laplacian H-eigenvalue is the spectral radius of the
signless Laplacian tensor \cite{cpz08,q12a,yy11}, since the signless
Laplacian tensor is a nonnegative tensor. There is sophisticated
theory for the spectral radius of a nonnegative tensor. Thus, the
corresponding theory for the largest signless Laplacian H-eigenvalue
is clear. On the other hand, the largest Laplacian H-eigenvalue is
more subtle. It can be seen that there are neat and simple
characterizations for the lower bound of the largest Laplacian
H-eigenvalue of an even-uniform hypergraph (Theorem \ref{thm-5}).
These are largely due to Lemma \ref{lem-5}. While, for odd-uniform
hypergraphs, the current theory is incomplete. This would be the
next topic to investigate.

{\bf Acknowledgement.} The authors are grateful to Prof. Jia-Yu Shao for his comments, and Prof. Xiaodong Zhang for Reference \cite{z07}.
\bibliographystyle{model6-names}

\end{document}